\documentclass[a4paper,10pt,openany]{article}
\usepackage{indentfirst}
\usepackage[greek,french,english]{babel}
\usepackage[iso-8859-7]{inputenc}

\usepackage{color,listings,titlesec}
\usepackage[headings]{fullpage}
\usepackage{amsfonts}
\usepackage{fancyhdr}
\usepackage{graphicx}
\usepackage{amsmath}
\newtheorem{Proposition}{Proposition}

\newtheorem{Theorem}[Proposition]{Theorem}

\newcommand{\g}{\mathfrak{g}}

\newcommand{\h}{\mathfrak{h}}

\begin{document}
\begin{center}

{\Large \textbf{Biquantization techniques for computing characters of differential operators on Lie groups.}}

\textbf{Panagiotis Batakidis}\footnote{Department of Mathematics, Aristotle University of Thessaloniki.}

\end{center}

\textbf{Abstract.} We compare the character of the algebra $(U(\g)/U(\g)\h_{\lambda})^{\h}$ as used by Fujiwara and Corwin-Greenleaf with the character produced from biquantization techniques applied in the Lie case by Cattaneo-Torossian. We prove that up to a smaller (specialization) algebra, these two characters are the same. An old example is also treated and it is proved that now we get more information about the question of when the symmetrization is an isomorphism of algebras.

\section{Introduction.}
 Let $G$ be a real nilpotent, connected and simply connected Lie group with $\g$ its Lie algebra, $\h$ a subalgebra of $\g$, $\lambda\in\h^{\ast}$ such that $\lambda([\h,\h])=0$. Then for $\;Y\in\h$, $\chi_{\lambda}:\;H\longrightarrow \mathbb{C}$ defined by $\chi_{\lambda}(\exp Y)=e^{i\lambda(Y)}$, is a unitary character of $H$ and we define the induced representation $\tau_{\lambda}:=Ind(G\uparrow H,\chi_{\lambda})$ with Hilbert space $\mathcal{H}_{\lambda}:=L^2(G,H,\lambda)$ the separable completion of $C^{\infty}_c(G,H,\chi_\lambda)$ with respect to the norm $||\phi||_2=\int_{G/H}|\phi(g)|^2d_{G/H}(g)$. The action of $G$ on $\phi\in L^2(G,H,\lambda)$ is translations by left: $\tau_{\lambda}(g)(\phi)(g')=\phi(g^{-1}g')$. These data correspond to a line bundle $\mathcal{L}_{\lambda}$ with base space $G/H$ and space of sections these functions $\phi$. Let $\mathcal{H}^{-\infty}_{\lambda}$ be the space of antilinear continuous forms on $\mathcal{H}_{\lambda}^{\infty}$, the later being the space of $C^{\infty}-$ vectors of $\mathcal{H}_{\lambda}$. The action of $U_{\mathbb{C}}(\g)$ on $\mathcal{H}^{-\infty}_{\lambda}$ will be denoted by $\mathrm{d}\tau_{\lambda}^{-\infty}$ and the action of $U_{\mathbb{C}}(\g)$ on $\mathcal{H}_{\lambda}^{\infty}$ is denoted respectively by $\mathrm{d}\tau_{\lambda}^{\infty}$.
 
 \textbf{Penney vectors.} Let $f\in\g^{\ast}$ such that $f|_{\h}=\lambda$ and $\mathfrak{b}$ a polarization with respect to $f$. We denote as $B$ its associated Lie group. Set $\h_{if}:=\langle H+if(H),\;H\in\h\rangle$. We shall denote as $U_{\mathbb{C}}(\g)\h_{-if}$ the ideal of $U_{\mathbb{C}}(\g)$ generated by $\h_{-if}$.  Let also $\mathrm{d}_{H,H\cap B}$ be a left-invariant measure on $H/H\cap B$. Let $\alpha_f$ be an element of $\mathcal{H}_f^{-\infty}$ defined for $\phi\in\mathcal{H}_f^{\infty}$, as 
 
 \begin{equation}\label{distribution}
<\alpha_f,\phi>=\int_{H/H\cap B}\overline{\phi(h)\chi_{\lambda}(h)}\mathrm{d}_{H/H\cap B}(h).
\end{equation}

The vector $\alpha_f$ is $H-$ semi-invariant (~\cite{B1}). Because of this invariance property of $\alpha_f$, the algebra $\left(U_{\mathbb{C}}(\g)/U_{\mathbb{C}}(\g)\h_{-if}\right)^{\h}$ is acting on $\alpha_f$. 

In an algebraic setting the lagrangian condition can be rewritten as
\begin{equation}\label{lagnow}
\exists \mathcal{O}\subset \lambda+\h^{\bot}\;\;\textlatin{a non-empty Zariski-open set, such that}\;\;\forall f\in\mathcal{O},\;\dim(\h\cdot f)=\frac{1}{2}\dim(\g\cdot f).
\end{equation} 
Recall that if the $H\cdot f$ orbits are lagrangian in the $G\cdot f$ orbits, then $(U_{\mathbb{C}}(\g)/U_{\mathbb{C}}(\g)\h_{-if})^{\h}$ is commutative (for this result see ~\cite{CG1} $\mathcal{x}$ 5, Theorem 5.4 and Corollary 5.5). 
 
\section{Construction of characters.}
One of our main objects is the reduction algebra  $H^0_{(\epsilon)}(\h_{\lambda}^{\bot},d^{(\epsilon)}_{\h_{\lambda}^{\bot},\mathfrak{q}})$. To briefly describe it we first need to describe the differential $d^{(\epsilon)}_{\h_{\lambda}^{\bot},\mathfrak{q}}:\;S(\mathfrak{q})[\epsilon]\longrightarrow S(\mathfrak{q})[\epsilon]\otimes\h^{\ast}$. This differential is defined as $d^{(\epsilon)}_{\h_{\lambda}^{\bot},\mathfrak{q}}:=\sum_{i=1}^{\infty}\epsilon^i d^{(i)}_{\h_{\lambda}^{\bot},\mathfrak{q}}$ where $d^{(i)}_{\h_{\lambda}^{\bot},\mathfrak{q}}:=\sum_{\Gamma\in \mathcal{B}_i\cup \mathcal{BW}_i}\overline{\omega}_{\Gamma}B_{\Gamma}$. Here $\Gamma$ stands for Kontsevich graphs (as in \cite{K}) that have to belong to $\mathcal{B}_i\cup \mathcal{BW}_i$, a special family of Kontsevich graphs (namely Bernoulli and Bernoulli attached to a wheel (see \cite{CT} for their description). The component $\overline{\omega}_{\Gamma}$ is a real coefficient depending on $\Gamma$ and $B_{\Gamma}$ is a differential operator depending also on $\Gamma$. For more details on the definitions and the formulas we refer to \cite{BAT} $\mathcal{x}$ 2.3.2, 3.2.1, or the note \cite{BAT2}.
The elements of this algebra are polynomials (on the formal deformation parameter $\epsilon$) $P_{(\epsilon)}$, which are solutions of the equation $d^{(\epsilon)}_{\h^{\bot}_{\lambda},\mathfrak{q}}(P_{(\epsilon)})=0$. The space of solutions is a vector space which we equip with the Cattaneo-Felder (associative) star-product $\ast_{CF,\epsilon}$ to take the reduction algebra $H^0_{(\epsilon)}(\h_{\lambda}^{\bot},d^{(\epsilon)}_{\h_{\lambda}^{\bot},\mathfrak{q}})$.

Choose a supplementary space $\mathfrak{q}$ for $\h$ in $\g$. Let $Y\in\g$, set $q(Y): = \det_{\g} \left(\frac{\sinh\frac{\mathrm{ad} Y}{2}}{\frac{adY}{2}}\right)$, and recall the symmetrization map $\beta:\;S(\g)\longrightarrow U(\g)$. Define $T_1,T_2$ to be the operators 

\[T_1:\;H^0_{(\epsilon)}(\g^{\ast},d^{(\epsilon)}_{\g^{\ast}})\longrightarrow H^0_{(\epsilon)}(\h_{\lambda}^{\bot},d^{(\epsilon)}_{\g^{\ast},\h_{\lambda}^{\bot},\mathfrak{q}}),\;\;\;\;F \mapsto F\ast_1 1,\] 
\[T_2:\;H^0_{(\epsilon)}(\h^{\bot}_{\lambda},d^{(\epsilon)}_{\h^{\bot}_{\lambda},\mathfrak{q}}) \longrightarrow H^0_{(\epsilon)}(\h_{\lambda}^{\bot},d^{(\epsilon)}_{\g^{\ast},\h_{\lambda}^{\bot},\mathfrak{q}}),\; \;\;\;G\mapsto 1 \ast_2 G.\]
that is the operators defining the Cattaneo-Felder bimodule structure on the biquantization diagramm of $\g^{\ast}$ and $\h_{\lambda}^{\bot}$. We denoted as $H^0_{(\epsilon)}(\h_{\lambda}^{\bot},d^{(\epsilon)}_{\g^{\ast},\h_{\lambda}^{\bot},\mathfrak{q}})$ the reduction space at the corner of this diagramm and as $\ast_1,\;\ast_2$ its left $H^0_{(\epsilon)}(\g^{\ast},d^{(\epsilon)}_{\g^{\ast}})$-module structure and its right $H^0_{(\epsilon)}(\h^{\bot}_{\lambda},d^{(\epsilon)}_{\h^{\bot}_{\lambda},\mathfrak{q}})$-module structure respectively. Using some simple facts we write them as $T_1:\left(S_{(\epsilon)}(\g),\ast_{DK}\right)\simeq \left(U_{(\epsilon)}(\g),\cdot\right)\longrightarrow S(\mathfrak{q})[\epsilon]$,  $\;T_2:\;H^0_{(\epsilon)}(\h^{\bot}_{\lambda},d^{(\epsilon)}_{\h^{\bot}_{\lambda},\mathfrak{q}})\longrightarrow S(\mathfrak{q})[\epsilon]$.

The PBW theorem holds for the deformed algebras $S_{(\epsilon)}(\g)=S(\mathfrak{q})[\epsilon]\oplus S_{(\epsilon)}(\g)\ast_{DK}\h_{\lambda}$ ($\ast_{DK}$ stands for the Duflo-Kontsevich star-product), and $U_{(\epsilon)}(\g)$ and there is a symmetrization map $\beta_{(\epsilon)}:\;S_{(\epsilon)}(\g)\longrightarrow~U_{(\epsilon)}(\g)$. We denote as $\overline{\beta}_{\mathfrak{q},(\epsilon)}:\;S(\mathfrak{q})[\epsilon]\longrightarrow U_{(\epsilon)}(\g)/U_{(\epsilon)}(\g)\h_{\lambda}$ the quotient of this symmetrization map with respect to the chosen $\mathfrak{q}$. We now write $\forall X\in\g,\;q_{(\epsilon)}(X):=q(\epsilon X)$ and note that using the isomorphism $(S_{(\epsilon)}(\g),\ast_{DK})\simeq (S_{(\epsilon)}(\g),\ast_{CF})\simeq(U_{(\epsilon)}(\g),\cdot)$,
 $U_{(\epsilon)}(\g)$ can be decomposed as $U_{(\epsilon)}(\g)=\bar{\beta}_{\mathfrak{q},(\epsilon)}\circ\partial_{q_{(\epsilon)}^{\frac{1}{2}}}(S(\mathfrak{q})[\epsilon])\oplus U_{(\epsilon)}(\g)\cdot\h_{\lambda}.$
Finally we will write $\overline{T}_1:=T_1|_{S(\mathfrak{q})[\epsilon]}$. It needs a small lemma to show that $\overline{T}_1$ is an isomorphism of vector spaces and we will denote abusively as $\overline{T}_1^{-1}$ its inverse $\overline{T}_1^{-1}:\;S(\mathfrak{q})[\epsilon]\longrightarrow S(\mathfrak{q})[\epsilon]\subset H^0_{(\epsilon)}(\g^{\ast},d^{(\epsilon)}_{\g^{\ast}})$. 

In \cite{BAT} $\mathcal{x}$ 3.4.2, Theorem 3.1 and the note \cite{BAT2} we proved that there is an explicit non-canonical isomorphism \[\overline{\beta}_{\mathfrak{q},(\epsilon)}\circ\partial_{q_{(\epsilon)}^{\frac{1}{2}}}\circ \overline{T}_1^{-1}T_2:\; H^0_{(\epsilon)}(\h^{\bot}_{\lambda},d^{(\epsilon)}_{\h^{\bot}_{\lambda},\mathfrak{q}})\stackrel{\simeq}{\longrightarrow} \left(U_{(\epsilon)}(\g)/U_{(\epsilon)}(\g)\h_{\lambda}\right)^{\h}.\]

 We shall use this fact to construct a family of characters by means of  \textlatin{Cattaneo- Felder - Torossian} techniques \cite{CF3}, \cite{CKTB}. More specifically,

\begin{Theorem}[see \cite{CT}]
Let $\g$ be a Lie algebra over $\mathbb{R}$, $\h\subset\g$,  $f\in\g^{\ast}$ s.t $\h$ is lagrangian with respect to $f$. Let  $\mathfrak{b}$ be a polarization of  $f$ and $\mathfrak{q}_{\mathfrak{b}}$ a transverse supplementary of $\h$.
The map
\[\gamma_{CT}:\;\left(U_{(\epsilon)}(\g)/U_{(\epsilon)}(\g)\h_{f}\right)^{\h}\longrightarrow\mathbb{R}[\epsilon]\]
\[u\mapsto \overline{T}_1^L\circ\bar{\beta}_{\mathfrak{q}_{\mathfrak{b}},(\epsilon)}^{-1}(u)(f)\]
is a character of  $(U_{(\epsilon)}(\g)/U_{(\epsilon)}(\g)\h_{\lambda})^{\h}$.
\end{Theorem}

Modifying the initial conditions we actually get something more useful: The previous theorem constructs only one character. However if the Lie group $G$ is nilpotent, under the general lagrangian condition ($\exists \mathcal{O}\subset -\lambda+\h^{\bot}$ s.t $\forall l\in\mathcal{O}$, the orbits $H\cdot l\subset G\cdot l$ are lagrangian submanifolds),  we can construct a character of  $(U_{(\epsilon)}(\g)/U_{(\epsilon)}(\g)\h_{\lambda})^{\h}$ for each such element $l$. Our goal here is to compare the character defined through the Penney eigendistribution in non-commutative harmonic analysis with that of deformation quantization.

Let $\mathcal{H}_f^{\infty}$ be the $C^{\infty}-$ vectors of the Hilbert space of the representation  $\tau_f=Ind(G,H,f)$ and $\alpha_f\in \mathcal{H}_f^{-\infty}$ the distribution defined for $\phi\in\mathcal{H}_f^{\infty}$ from the formula $<\alpha_f,\phi>=\int_{H/H\cap B}\overline{\phi(h)\chi_{\lambda}(h)}\mathrm{d}_{H/H\cap B}(h)$.

\begin{Theorem}[\cite{FUJI2}]
Let $\g$ be a Lie algebra $(\dim(\g)<\infty)$, $\h\subset\g$, $\lambda$ a character of $\h$.  Suppose that generically the representation $\tau_{\lambda}=Ind(G,H,\lambda)$ has finite multiplicities in her spectral decomposition. Then for $l\in \lambda+\h^{\bot}$ and $A\in(U_{\mathbb{C}}(\g)/U_{\mathbb{C}}(\g)\h_{il})^{\h}$, the action $\mathrm{d}\tau_l(\overline{A})(\alpha_l)$ is a multiple of  $\alpha_l$, Thus there exists a character $\lambda_l:\;(U_{\mathbb{C}}(\g)/U_{\mathbb{C}}(\g)\h_{il})^{\h}\longrightarrow \mathbb{C}$ defined from the relation $\mathrm{d}\tau_l(\overline{A})(\alpha_l)=\overline{\lambda_l(A)}\alpha_l$.
\end{Theorem}

Before we proceed, we need to define two specialization algebras. First we set $H^0_{(\epsilon=1)}(\h^{\bot}_{\lambda},d^{(\epsilon=1)}_{\h^{\bot}_{\lambda},\mathfrak{q}}):=H^0_{(\epsilon)}(\h^{\bot}_{\lambda},d^{(\epsilon)}_{\h^{\bot}_{\lambda},\mathfrak{q}})/\langle\epsilon-1\rangle$ to be the specialization algebra of the reduction algebra $H^0_{(\epsilon)}(\h^{\bot}_{\lambda},d^{(\epsilon)}_{\h^{\bot}_{\lambda},\mathfrak{q}})$.

Consider a supplementary variable $T$ such that $[T,\g]=0$ and set $\g_T=\g\oplus <T>$ and $\h_T=\h\oplus <T>$ such that $\dim(\g_T)=\dim(\g)+1$. Set also $U(\g_T)$ to be the universal enveloping algebra of $\g_T$ and $U(\g_T)\h^T_{\lambda}$ to be the ideal of $U(\g_T)$ generated by  $\h_{\lambda}^T=<H+T\lambda(H),\;H\in\h>$. Let $H$ be the associated Lie group of $\h$ and consider the unitary character $\chi_{\lambda}:\;H\longrightarrow \mathbb{C}$ defined by the formula for $Y\in\h$, $\chi_{\lambda}(\exp Y)= \exp(i\lambda(Y))$. Denote as $C^{\infty}(G,H,\chi_\lambda)$ the vector space of complex smooth functions $\theta$ on $G$ that satisfy the property $\forall h\in H, \forall g\in G,\;\;\theta (gh)=\chi^{-1}_{\lambda}(h)\theta(g)$. We denote as $\mathbb{D}(\g,\h,\lambda)$ the algebra of linear differential operators, that leave the space $C^{\infty}(G,H,\chi_\lambda)$ invariant and commute with the left translation on $G$. 

Recall that from a theorem of Koornwider we have $(U(\g)/U(\g)\h_{\lambda})^{\h}\simeq\mathbb{D}(\g,\h,\lambda)$. Thus setting $\mathbb{D}_T(\g,\h,\lambda):= \mathbb{D}(\g_T,\h_T,\lambda)$ we can also write $(U(\g_T)/U(\g_T)\h^T_{\lambda})^{\h_T}\simeq\mathbb{D}_T(\g,\h,\lambda)$. Finally we define our second specialization algebra $\mathbb{D}_{(T=1)}(\g,\h,\lambda):=(U(\g_T)/U(\g_T)\h^T_{\lambda})^{\h_T}/\langle T-1\rangle$. In \cite{BAT} $\mathcal{x}$ 3.5.3 Theorem 3.5 it is proved as the outcome of series of other results that $\mathbb{D}_{(T=1)}(\g,\h,\lambda)\simeq H^0_{(\epsilon=1)}(\h^{\bot}_{\lambda},d^{(\epsilon=1)}_{\h^{\bot}_{\lambda},\mathfrak{q}})$. This result can also be found at \cite{BAT2}.

In order to proceed to the character comparison, it is necessary that the theorems of harmonic analysis and deformation quantization refer to the same field. So we need a real character since the whole Kontsevich construction which we mentioned is over $\mathbb{R}$:

\begin{Theorem}[\cite{BAT}, $\mathcal{x}$ 4.4.2, Theorem 4.3]
Let $\g,\h,\lambda$ as before and suppose that the $H-$ orbits are lagrangian in the affine space $\lambda+\h^{\bot}$. Then for a regular $f\in \lambda+\h^{\bot}$ and s.t  $\dim(\h\cdot f)=\frac{1}{2}\dim(\g\cdot f)$, and $A\in D_{(T=1)}(\g,\h,\lambda)$, the action $\mathrm{d}\tau_f(A)(\alpha(f))$is a multiple of $\alpha(f)$, and so there is defined a character $\lambda_{(T=1)}^f:\;D_{(T=1)}(\g,\h,\lambda)\longrightarrow \mathbb{R}$ such that $\mathrm{d}\tau_f(A)(\alpha(f))=\lambda_{(T=1)}^f(A)\alpha(f).$
\end{Theorem}

\textbf{Proof.}  The idea is to follow the line of proof of Fujiwara proving theorem 2. This is done by double induction on $\dim(\g)$ and $\dim(\h)$ and works fine up to the case $\h\subset\g_0$, $\g_0$ being a codimension one ideal of $\g$ (constructed in a standard way using the reduction triplet in the sense of Dixmier). In this case, the condition of Corwin-Greenleaf (see (1) of equations (2.7) in ~\cite{FLMM}) holds for the character $it\lambda$ and we have 
\begin{equation}\label{ratiot}
\left(U_{\mathbb{C}}(\g)/U_{\mathbb{C}}(\g)\h_{it\lambda}\right)^{\h}=\left(U_{\mathbb{C}}(\g_0)/U_{\mathbb{C}}(\g_0)\h_{it\lambda}\right)^{\h}.
\end{equation}

This equation depends rationally on $it$, $t\in\mathbb{R}^{\ast}$. So if (\ref{ratiot}) holds for $it$, $t\in\mathbb{R}^{\ast}$, it holds also for $t$ in a Zariski-open subset of $\mathbb{R}$ and we write $\left(U_{\mathbb{C}}(\g)/U_{\mathbb{C}}(\g)\h_{t\lambda}\right)^{\h}=\left(U_{\mathbb{C}}(\g_0)/U_{\mathbb{C}}(\g_0)\h_{t\lambda}\right)^{\h}$, but we can't conclude that a similar equation holds for the algebra $\left(U_{\mathbb{C}}(\g)/U_{\mathbb{C}}(\g)\h_{\lambda}\right)^{\h}$. This is the difference with respect to the proof of Thm. 2 which was about a unitary character. To overcome this setback, we can use polynomial families $t\mapsto u_t\in  \left(U_{\mathbb{C}}(\g_0)/U_{\mathbb{C}}(\g_0)\h_{t\lambda}\right)^{\h}$. This will allow us to continue the argument and it explains at the same time why the character at the theorem's statement is defined for $D_{(T=1)}(\g,\h,\lambda)$ which is elsewhere (\cite{BAT} $\mathcal{x}$ 3.5.3, Cor. 3.2) shown to correspond to elements who are the value at $t=1$ of polynomial families $t\mapsto u_t\in \left(U(\g)/U(\g)\h_{t\lambda}\right)^{\h}$. The proof then continues by carefully applying the induction arguments on $\dim(\g)$. $\diamond$

Thus the real character that we construct has a price: The corresponding theorem for the Penney distribution now holds for a smaller algebra: Here becomes also clear the use of the specialization algebra $D_{(T=1)}(\g,\h,\lambda)\simeq H^0_{(\epsilon=1)}(\h_{\lambda}^{\bot},d^{(\epsilon=1)}_{\h_{\lambda}^{\bot},\mathfrak{q}})$  introduced in \cite{BAT}. This algebra might be the appropriate object of study when it comes to the Duflo and Corwin-Greenleaf conjectures. 

\section{Comparison of characters and example.}
Let $\mathfrak{i}_{(\epsilon=1)}:\;H^0_{(\epsilon=1)}(\h_{\lambda}^{\bot},d^{(\epsilon=1)}_{\h_{\lambda}^{\bot},\mathfrak{q}})\hookrightarrow(U(\g)/U(\g)\h_{\lambda})^{\h}$ the injective map coming from the fact that $H^0_{(\epsilon=1)}(\h_{\lambda}^{\bot},d^{(\epsilon=1)}_{\h_{\lambda}^{\bot},\mathfrak{q}})\hookrightarrow H^0(\h_{\lambda}^{\bot},d_{\h_{\lambda}^{\bot},\mathfrak{q}})$ and $H^0(\h_{\lambda}^{\bot},d_{\h_{\lambda}^{\bot},\mathfrak{q}})\hookrightarrow (U(\g)/U(\g)\h_{\lambda})^{\h}$.

\begin{Theorem}[\cite{BAT}, $\mathcal{x}$ 4.4.4, Theorem 4.4]
Let $\g$ a nilpotent Lie algebra ($\dim(\g)<\infty$), $\h\subset \g$ a subalgebra, $\lambda$ a character of  $\h$. Let $P\in  H^0_{(\epsilon=1)}(\h_{\lambda}^{\bot},d^{(\epsilon=1)}_{\h_{\lambda}^{\bot},\mathfrak{q}})$ and $u\in D_{(T=1)}(\g,\h,\lambda)$ s.t $u= \mathfrak{i}_{(\epsilon=1)}(P)$. Then for a generic $f\in\lambda+\h^{\bot}$ there is a pair $(\mathfrak{b}_f, \mathfrak{q}_f)$ satysfying $T_1\circ\overline{\beta}^{-1}_{\mathfrak{q}_f,(\epsilon)}(P)|_{\epsilon=1}(-f)=\lambda^f_{(T=1)}(u).$
\end{Theorem}

\textbf{Proof.} This is done again by a long double induction on $\dim(\g)$ and $\dim(\h)$ confirming that in every step, we compute in the same subspaces for $D_{(T=1)}(\g,\h,\lambda)$ and $H^0_{(\epsilon=1)}(\h_{\lambda}^{\bot},d^{(\epsilon=1)}_{\h_{\lambda}^{\bot},\mathfrak{q}})$ and that the computations match, giving the same character. $\diamond$

\textbf{Example.} We end the present note with an example that reveals the power of this approach (the fully detailed and computed example is at  $\mathcal{x} 5.5$ of \cite{BAT}):
Let $\mathfrak{g}$ be the nilpotent Lie algebra generated by ${X,U,V,E,Z}$ with relations $[U,V]=E, [X,U]=V, [X,V]=Z$,
$\mathfrak{h}=\mathbb{R}X\oplus \mathbb{R}E$ and $\lambda=E^{\ast}$. 
For a $u\in (U(\g)/U(\g)\h_{\lambda})^{\h}$ and with the right choices (transversal condition) of $\mathfrak{q}_l,\mathfrak{q}$ we have $\beta_{\mathfrak{q}_l}^{-1}(u)=e^{[\frac{1}{12l(Z)}(1-\frac{Z}{2l(Z)})\partial_U^3]}\beta_{\mathfrak{q}}^{-1}(u),$ where $\beta_{\mathfrak{q}}$ is the quotient symmetrization map. If $v$ is a polynomial of $(U(\g)/U(\g)\h_{\lambda})^{\h}$ then $\beta_{\mathfrak{q}}^{-1}(v)(l)=e^{-\frac{1}{24l(Z)}\partial_U^3}\beta^{-1}_{\mathfrak{q}_l}(v)(l)$.

The map $\gamma_{CT}:\;v\mapsto \left(e^{[\frac{1}{12l(Z)}(1-\frac{Z}{2l(Z)})\partial_U^3]}\beta_{\mathfrak{q}}^{-1}(v)\right)(l)$ is a character of the algebra of differential operators $(U(\g)/U(\g)\h_{\lambda})^{\h}$. The important point here is the term in the exponential. This example was long before treated as a counterexample to the idea that the symmterization map $\beta$ was an algebra isomorphism in this case. Despite the example thought it was not possible to compute its exact formula of the isomorphism. The relation $\beta_{\mathfrak{q}_l}^{-1}(u)=e^{[\frac{1}{12l(Z)}(1-\frac{Z}{2l(Z)})\partial_U^3]}\beta_{\mathfrak{q}}^{-1}(u)$ reveals the problem and computes in this case the extra term of third degree with rational coefficients which is only computed using the deformation quantization techniques.

\textbf{Acknowledgement.} The author would like to gratefully thank Charles Torossian for his support, inspiration and guiding excellence during his PhD thesis at University Paris 7.

\end{document}